\documentclass[12pt]{amsart}
\usepackage{amsmath,amssymb,latexsym}
\def\P{\mathcal P_\alpha}
\def\SU{\mathrm{SU}}
\def\Pcl{\Sigma\mathrm{cl}}
\def\PPcl{\Sigma'\mathrm{cl}}
\def\Qcl{\Sigma_1\mathrm{cl}}
\def\Acl{\Sigma_\alpha\mathrm{cl}}
\def\dcl{\mathrm{dcl}}
\def\acl{\mathrm{acl}}
\def\bdd{\mathrm{bdd}}
\def\cb{\mathrm{Cb}}

\def\M{\mathfrak M}

\def\tp{\mathrm{tp}}
\def\lstp{\mathrm{Lstp}}

\def\domeq{\mathbin{\underline\square}}
\def\Ind#1#2{#1\setbox0=\hbox{$#1x$}\kern\wd0\hbox to 0pt{\hss$#1\mid$\hss}
\lower.9\ht0\hbox to 0pt{\hss$#1\smile$\hss}\kern\wd0}
\def\ind{\mathop{\mathpalette\Ind{}}}
\def\Notind#1#2{#1\setbox0=\hbox{$#1x$}\kern\wd0\hbox to 0pt{\mathchardef
\nn="3236\hss$#1\nn$\kern1.4\wd0\hss}\hbox to 0pt{\hss$#1\mid$\hss}\lower.9\ht0
\hbox to 0pt{\hss$#1\smile$\hss}\kern\wd0}
\def\nind{\mathop{\mathpalette\Notind{}}}

\theoremstyle{plain}
\newtheorem{theorem}{Theorem}[section]
\newtheorem{proposition}[theorem]{Proposition}
\newtheorem{fact}[theorem]{Fact}
\newtheorem{lemma}[theorem]{Lemma}
\newtheorem{corollary}[theorem]{Corollary}
\newtheorem*{claim}{Claim}

\theoremstyle{definition}
\newtheorem{definition}[theorem]{Definition}
\newtheorem{remark}[theorem]{Remark}
\newtheorem{question}[theorem]{Question}
\newtheorem{problem}[theorem]{Problem}
\newtheorem*{expl}{Example}
\newtheorem*{conj}{Conjecture}

\def\bsp{\begin{expl}}
\def\ebsp{\end{expl}}
\def\verm{\begin{conj}}
\def\everm{\end{conj}}
\def\beh{\begin{claim}}
\def\ebeh{\end{claim}}
\def\defn{\begin{definition}}
\def\edefn{\end{definition}}
\def\satz{\begin{theorem}}
\def\esatz{\end{theorem}}
\def\tats{\begin{fact}}
\def\etats{\end{fact}}
\def\kor{\begin{corollary}}
\def\ekor{\end{corollary}}
\def\lmm{\begin{lemma}}
\def\elmm{\end{lemma}}
\def\frag{\begin{question}}
\def\efrag{\end{question}}
\def\bem{\begin{remark}}
\def\ebem{\end{remark}}
\def\prob{\begin{problem}}
\def\eprob{\end{problem}}
\def\bew{\par\noindent{\em Proof: }}
\def\bewbeh{\par\noindent{\em Proof of Claim: }}
\def\satzli{\begin{proposition}}
\def\esatzli{\end{proposition}}

\begin{document}
\title{Ample thoughts}
\author{Daniel Palac\'\i n and Frank O. Wagner}
\thanks{The first author was partially supported
by research project MTM 2008-01545 of the Spanish government and
research project 2009SGR 00187 of the Catalan government. He also would like to thank Amador Mart\' in Pizarro for interesting discussions around the Canonical Base Property, and suggesting a simplification in the proof of Theorem \ref{domequ}.\\
The second author was partially supported by ANR-09-BLAN-0047 Modig}
\address{Universitat de Barcelona; Departament de L\`ogica, Hist\`oria i Filosofia de la Ci\`encia, Montalegre 6, 08001 Barcelona, Spain}
\address{{\it Current address:} Universit\'e de Lyon; CNRS; Universit\'e Lyon 1; Institut Camille Jordan UMR5208, 43 bd du 11 novembre 1918, 69622 Villeurbanne Cedex, France}
\address{Universit\'e de Lyon; CNRS; Universit\'e Lyon 1; Institut Camille Jordan UMR5208, 43 bd du 11 novembre 1918, 69622 Villeurbanne Cedex, France}
\email{palacin@math.univ-lyon1.fr}
\email{wagner@math.univ-lyon1.fr}
\keywords{stable; simple; one-based; CM-trivial; $n$-ample; internal; analysable; closure; level; flat; ultraflat; canonical base property}
\subjclass[2000]{03C45}
\date{4 September 2012}
\begin{abstract}Non-$n$-ampleness as defined by Pillay \cite{pi00} and Evans \cite{ev03} is preserved under analysability. Generalizing this to a more general notion of $\Sigma$-ampleness, this gives an immediate proof for all simple theories of a weakened version of the Canonical Base Property (CBP) proven by Chatzidakis \cite{zoe} for types of finite SU-rank. This is then applied to the special case of groups.\end{abstract}
\maketitle

\section{Introduction}
Recall that a partial type $\pi$ over a set $A$ in a simple theory is {\em
one-based} if for any tuple $\bar a$ of realizations of $\pi$ and any
$B\supseteq A$ the canonical base $\cb(\bar a/B)$ is contained in the bounded closure
$\bdd(\bar aA)$. In other words, forking dependence is either trivial or behaves as in modules: Any two sets are independent over the intersection of their bounded closures. One-basedness implies that the forking geometry is
particularly well-behaved; for instance one-based groups are
bounded-by-abelian-by-bounded. The principal result in \cite{wa04} is that one-basedness is preserved under analyses (i.e.\ iterative approximations by some other types): a type analysable in one-based types is itself one-based. This generalized earlier results of Hrushovski \cite{udi} and Chatzidakis \cite{zoe}. One-basedness is the first level in a hierarchy of possible geometric behaviour of forking independence first defined by Pillay \cite{pi00} and slightly modified by Evans \cite{ev03}, $n$-ampleness, modelled on the behavior of flags in $n$-space. Not $1$-ample means one-based; not $2$-ample is equivalent to a notion previously introduced by Hrushovski \cite{udi93}, CM-triviality. Fields are $n$-ample for all $n<\omega$, as is the non-abelian free group \cite{oht}. In \cite{pi00} Pillay defines $n$-ampleness locally for a single type and shows that a superstable theory of finite Lascar rank is non $n$-ample if and only if all its types of rank $1$ are; his proof implies that in such a theory, a type analysable in non $n$-ample types is again non $n$-ample.

We shall give a definition of $n$-ampleness for invariant families of partial types, and generalize Pillay's result to arbitrary simple theories. Note that for $n=1$ this gives an alternative proof of the main result in \cite{wa04}. Since for types of infinite rank the algebraic (bounded) closure used in the definition is not necessarily appropriate (for a regular type $p$ one might, for instance, replace it by $p$-closure), we also generalize the notion to $\Sigma$-closure for some $\emptyset$-invariant collection of partial types (thought of as small), giving rise to the notion of $n$-$\Sigma$-ample. This may for instance be applied to consider ampleness modulo types of finite $\SU$-rank, or modulo
supersimple types. Readers not interested in this additional generality are invited to simply replace $\Sigma$-closure by bounded closure. However, this will only marginally shorten the proofs. As an immediate Corollary of the more general version, we shall derive a weakened version of the Canonical Base Property CBP \cite{PZ03} shown by Chatzidakis \cite{zoe}, where analysability replaces internality in the definition. We also give a version appropriate for supersimple theories. Finally, we deduce that in a simple theory with enough regular types, a hyperdefinable group modulo its approximate centre is analysable in the family of non one-based regular types; the group modulo a normal nilpotent subgroup is almost internal in that family. This can be thought of as a general version of the properties of one-based groups mentioned above.

Our notation is standard and follows \cite{wa00}.
Throughout the paper, the ambient theory will be simple, and we shall
be working in $\M^{heq}$, where $\M$ is a sufficiently saturated model
of the ambient theory. Thus tuples are tuples of hyperimaginaries, and
$\dcl=\dcl^{heq}$.

\section{Internality and analysability}

For the rest of the paper $\Sigma$ will be an $\emptyset$-invariant family of partial types. Recall first the definitions of internality, analysability, foreignness and orthogonality.
\defn Let $\pi$ be a partial type over $A$. Then $\pi$ is\begin{itemize}
\item ({\em almost}) {\em $\Sigma$-internal} if for every realization $a$
of $\pi$ there is $B\ind_Aa$ and a tuple $\bar b$ of realizations of types in $\Sigma$
based on $B$, such that $a\in\dcl(B\bar b)$ (or $a\in\bdd(B\bar b)$,
respectively).
\item {\em $\Sigma$-analysable} if for any realization $a$ of $\pi$ there are
$(a_i:i<\alpha)\in\dcl(Aa)$ such that $\tp(a_i/A,a_j:j<i)$ is
$\Sigma$-internal for all $i<\alpha$, and
$a\in\bdd(A,a_i:i<\alpha)$.\end{itemize}
A type $\tp(a/A)$ is {\em foreign} to $\Sigma$ if $a\ind_{AB}\bar b$
for all $B\ind_Aa$ and $\bar b$ realizing types in $\Sigma$ over
$B$. \newline Finally, $p\in S(A)$ is {\em orthogonal} to $q\in S(B)$ if for all $C\supseteq AB$, $a\models p$, and $b\models q$ with $a\ind_A C$ and $b\ind_B C$ we have $a\ind_C b$. \edefn

So $p$ is foreign to $\Sigma$ if $p$ is orthogonal to all completions of partial types in $\Sigma$, over all possible parameter sets.

The following lemmas and their corollaries are folklore, but we add some precision about non-orthogonality.
\lmm\label{orthogonal} Suppose $a\ind b$ and $a\nind_b c$. Let $(b_i:i<\omega)$ be an indiscernible sequence in $\tp(b)$
and put $p_b=\tp(c/b)$. Then $p_{b_i}$ is non-orthogonal to $p_{b_j}$ for all $i,j<\omega$.\elmm
\bew We prolong the sequence to have length $\alpha$. As $a\ind b$ and $(b_i:i<\alpha)$ is indiscernible, by \cite[Theorem 2.5.4]{wa00} we may assume $ab\equiv ab_i$ for all $i<\alpha$ and
$a\ind(b_i:i<\alpha)$. Let $B=(b_i:i<\omega)$, so $(b_i:\omega\le i<\alpha)$ is independent over $B$ and $a\ind B$. Choose $c_i$ with $b_ic_i\equiv_a bc$ and
$$c_i\ind_{ab_i}(b_j:j<\alpha)$$
for all $\omega\le i<\alpha$. Then $ab_ic_i\ind_B(b_j:j\not=i)$ for all $\omega\le i<\alpha$. By indiscernability, if $p_{b_i}$ were orthogonal to $p_{b_j}$ for some $i\not=j$, then they would be orthogonal for all $i\not=j$. As $c_i\ind_{b_i}(b_j:j\not=i)$, the sequence $(b_ic_i:\omega\le i<\alpha)$ would be independent over $B$. However, $a\nind_B b_ic_i$ for all $\omega\le i<\alpha$, contradicting the boundedness of weight of $\tp(a/B)$.\qed

\lmm\label{cbase} Suppose $a\ind b$ and $a'=\cb(bc/a)$. Let $\mathcal P$ be the family of $\bdd(\emptyset)$-conjugates of $\tp(c/b)$ non-orthogonal to $\tp(c/b)$. Then $a'\in\bdd(a)$ is $\mathcal P$-internal and $\bdd(ab)\cap\bdd(bc)\subseteq\bdd(a'b)$.\elmm
\bew If $a\ind bc$ then $a'\in\bdd(\emptyset)$ and $\bdd(ab)\cap\bdd(bc)=\bdd(b)$, so there is nothing to show. Assume $a\nind_bc$. Clearly $a'\in\bdd(a)$; as $bc\ind_{a'}a$ we get $c\ind_{a'b}a$ and hence $\bdd(ab)\cap\bdd(bc)\subseteq\bdd(a'b)$. Let $(b_ic_i:i<\omega)$ be a Morley sequence in $\lstp(bc/a)$ with $b_0c_0=bc$. Then $a'\in\dcl(b_ic_i:i<\omega)$; since $b\ind a$ we get $(b_i:i<\omega)\ind a$, whence $(b_i:i<\omega)\ind a'$. So $a'$ is internal in $\{\tp(c_i/b_i):i<\omega\}$. Finally, $\tp(c_i/b_i)$ is non-orthogonal to $\tp(c/b)$ for all $i<\omega$ by Lemma \ref{orthogonal}.\qed

\kor\label{cbase1} If $a\ind b$ and $\tp(c/b)$ is (almost) $\Sigma$-internal, then $\cb(bc/a)$ is (almost) $\Sigma$-internal. The same statement holds with {\em analysable} instead of {\em internal}.\ekor
\bew Let $d\ind_bc$ and $\bar e$ realize partial types in $\Sigma$ over $bd$ such that $c\in\dcl(bd\bar e)$ (or $c\in\bdd(bd\bar e)$, respectively). We may take $d\bar e\ind_{bc}a$. Then $d\ind_bac$, whence $a\ind bd$. So $\cb(bd\bar e/a)$ is $\Sigma$-internal by Lemma \ref{cbase}.
But $a\ind_{bc}d\bar e$ and $c\in\dcl(bd\bar e)$ implies $\cb(bc/a)\in\dcl(\cb(bd\bar e/a))$; similarly $c\in\bdd(bd\bar e)$ implies $\cb(bc/a)\in\bdd(\cb(bd\bar e/a))$.

The proof for $\Sigma$-analysability is analogous.\qed
\defn Two partial types $\pi_1$ and $\pi_2$ are {\em perpendicular}, denoted $\pi_1\ \underline\perp\ \pi_2$, if for any set $A$ containing their domains and any tuple $\bar a_i\models\pi_i$ for $i=1,2$ we have $\bar a_1\ind_A\bar a_2$.\edefn
For instance, orthogonal types of rank 1 are perpendicular.
\kor Suppose $a\ind b$, and $a_0\in\bdd(ab)$ is (almost) $\Pi$-internal over $b$ for some $b$-invariant family $\Pi$ of partial types. Let $\Pi'$ be the family of $\bdd(\emptyset)$-conjugates $\pi'$ of partial types $\pi\in\Pi$ with $\pi'\not\!\!\underline\perp\ \pi$. Then there is (almost) $\Pi'$-internal $a_1\in\bdd(a)$ with $a_0\in\bdd(a_1b)$. The same statement holds with {\em analysable} instead of {\em internal}.\ekor
\bew If $\tp(a_0/b)$ is $\Pi$-internal, there is $c\ind_ba_0$ and $\bar e$ realizing partial types in $\Pi$ over $bc$ such that $a_0\in\dcl(bc\bar e)$; we choose them with $c\bar e\ind_{ba_0}a$. So $c\ind_ba$, whence $a\ind bc$. Furthermore, we may assume that $e\nind_{bc}a$ for all $e\in\bar e$, since otherwise $ec\ind_ba_0$ and we may just include $e$ in $c$, reducing the length of $\bar e$. Now $a_0\in\bdd(abc)\cap\bdd(bc\bar e)$, so by Lemmas \ref{orthogonal} and \ref{cbase} there is $\Pi'$-internal $a_1\in\bdd(a)$ with $a_0\in\bdd(bca_1)$. Since $a\ind_bc$ implies $a_0\ind_{a_1b}c$, we get $a_0\in\bdd(a_1b)$.

For the almost internal case, we replace definable by bounded closure; for the analysability statement we iterate, adding $a_1$ to the parameters.\qed

To finish this section, a decomposition lemma for almost internality.
\lmm\label{intblock} Let $\Sigma=\bigcup_{i<\alpha}\Sigma_i$, where $(\Sigma_i:i<\alpha)$ is a collection of pairwise perpendicular $\emptyset$-invariant families of partial types. If $\tp(a/A)$ is almost $\Sigma$-internal, then there are $(a_i:i<\alpha)$ interbounded over $A$ with $a$ such that $\tp(a_i/A)$ is $\Sigma_i$-internal for $i<\alpha$.\elmm
Clearly, if $a$ is a finite imaginary tuple, we only need finitely many~$a_i$.
\bew By assumption there is $B\ind_A a$ and some tuples $(b_i:i<\alpha)$ such that $b_i$ realizes partial types in $\Sigma_i$ over $B$, with $a\in\bdd(B,b_i:i<\alpha)$. Let $a_i=\cb(Bb_i/Aa)$. Then $a_i\in\bdd(Aa)$ and $\tp(a_i/A)$ is $\Sigma_i$-internal by Corollary \ref{cbase1}.

Put $\bar a=(a_i:i<\alpha)$. Then $a\ind_{Aa_i}Bb_i$ implies $a\ind_{B\bar a}b_i$; since $b_i\ind_{Ba}(b_j:j\not=i)$ by perpendicularity we obtain $b_i\ind_{B\bar a}(a,b_j:j\not=i)$ for all $i<\alpha$. Hence $(a,b_i:i<\alpha)$ is independent over $B\bar a$, and in particular
$$a\ind_{B\bar a}(b_i:i<\alpha).$$
Since $a\in\bdd(B,b_i:i<\alpha)$ we get $a\in\bdd(B\bar a)$; as $a\ind_AB$ implies $a\ind_{A\bar a}B$ we obtain $a\in\bdd(A\bar a)$.\qed

\section{$\Sigma$-closure, $\Sigma_1$-closure and a theory of levels}

In his proof of Vaught's conjecture for superstable theories of finite rank \cite{bue08}, Buechler defines the first level $\ell_1(a)$ of an element $a$ of finite Lascar rank as the set of all $b\in\acl^{eq}(a)$ internal in the family of all types of Lascar rank one; higher levels are defined inductively by $\ell_{n+1}(a)=\ell_1(a/\ell_n(a))$. The notion  has been studied by Prerna Bihani Juhlin in her thesis \cite{prerna-thesis} in connection with a reformulation of the canonical base property. We shall generalise the notion to arbitrary simple theories.
\defn For an ordinal $\alpha$ the {\em $\alpha$-th $\Sigma$-level} of $a$ over $A$ is defined inductively by $\ell_0^\Sigma(a/A)=\bdd(A)$, and for $\alpha>0$
$$\ell_\alpha^\Sigma(a/A)=\{b\in\bdd(aA):\tp(b/\bigcup_{\beta<\alpha}\ell_\beta(a/A))\text{ is almost $\Sigma$-internal}\}.$$
Finally, we shall write $\ell_\infty^\Sigma(a/A)$ for the set of all hyperimaginaries $b\in\bdd(aA)$ such that $\tp(b/A)$ is $\Sigma$-analysable. \edefn
\bem Clearly, $\tp(a/A)$ is $\Sigma$-analysable if and only if $\ell_\infty^\Sigma(a/A)=\bdd(aA)$ if and only if $\ell_\alpha^\Sigma(a/A)=\bdd(aA)$ for some ordinal $\alpha$, and the minimal such $\alpha$ is the minimal length of a $\Sigma$-analysis of $a$ over $A$.\ebem

\lmm\label{level1-ind} If $a\ind b$, then $\ell_\alpha^\Sigma(ab)=\bdd(\ell_\alpha^\Sigma(a),\ell_\alpha^\Sigma(b))$ for any $\alpha$. \elmm
\bew Let $c=\ell_\alpha^\Sigma(ab)$. Clearly, $\ell_\alpha^\Sigma(a)\ell_\alpha^\Sigma(b)\subseteq c$. Conversely, put $a_0=\cb(bc/a)$. Then $\tp(a_0)$ is internal in the family of $\bdd(\emptyset)$-conjugates of $\tp(c/b)$ by Corollary \ref{cbase1}; since even $\tp(c)$ is $\Sigma$-analysable in $\alpha$ steps, so is $\tp(a_0)$. Thus $a_0\subseteq\ell_\alpha^\Sigma(a)$. Now $bc\ind_{a_0} a$ implies $$c\ind_{\ell_\alpha^\Sigma(a)b} a,$$ whence $c\subseteq\bdd(\ell_\alpha^\Sigma(a),b)$. By symmetry, $c\subseteq\bdd(\ell_\alpha^\Sigma(b),a)$, that is,
$$\ell_\alpha^\Sigma(ab)\subseteq \bdd(\ell_\alpha^\Sigma(a),b)\cap\bdd(\ell_\alpha^\Sigma(b),a).$$
On the other hand, $a\ind b$ yields $a\ind_{\ell_\alpha^\Sigma(a)\ell_\alpha^\Sigma(b)} b.$ Thus, $$\bdd(\ell_\alpha^\Sigma(a),b)\cap\bdd(\ell_\alpha^\Sigma(b),a)=\bdd(\ell_\alpha^\Sigma(a),\ell_\alpha^\Sigma(b)),$$ whence the result. \qed

\kor\label{cor-level1-ind} If $(a_i:i\in I)$ is an $\emptyset$-independent sequence, then $\ell_\alpha^\Sigma(a_i:i\in I)=\bdd(\ell_\alpha^\Sigma(a_i):i\in I)$. \ekor
\bew Let $c=\ell_\alpha^\Sigma(a_i:i\in I)$ and set $b_J=\cb(c/a_i:i\in J)$ for each finite $J\subseteq I$. Note for each finite $J\subseteq I$ that $\tp(b_J)$ is $\Sigma$-analysable in $\alpha$ steps. Thus $b_J\subseteq\ell_\alpha^\Sigma(a_i:i\in J)$. On the other hand,
$$\ell_\alpha^\Sigma(a_i:i\in J)=\bdd(\ell_\alpha^\Sigma(a_i):i\in J)$$ 
by Lemma \ref{level1-ind} and induction, since $J\subseteq I$ is finite. Therefore 
$$c\ind_{(\ell_\alpha^\Sigma(a_i):i\in I)} (a_i:i\in I)$$
by the finite character of forking, whence $c\subseteq\bdd(\ell_\alpha^\Sigma(a_i):i\in I)$.
\qed

We shall see that the first level governs domination-equivalence.
\defn An element $a$ {\em $\Sigma$-dominates} an element $b$ over $A$, denoted $a\unrhd_A^\Sigma b$, if for all $c$ such that $\tp(c/A)$ is $\Sigma$-analysable, $a\ind_Ac$ implies $b\ind_Ac$. Two elements $a$ and $b$ are {\em $\Sigma$-domination-equivalent} over $A$, denoted $a\domeq_A^\Sigma b$, if $a\unrhd_A^\Sigma b$ and $b\unrhd_A^\Sigma a$. If $\Sigma$ is the set of all types, it is omitted.\edefn
The following generalizes a theorem by Buechler \cite[Proposition 3.1]{bue08} for finite Lascar rank.
\satz\label{domequ} Let $\Sigma'$ be an $\emptyset$-invariant family of partial types.\begin{enumerate}
\item $a$ and $\ell_1^\Sigma(a/A)$ are $\Sigma$-domination-equivalent over $A$.
\item If $\tp(a/A)$ is $\Sigma$-analysable, then $a$ and $\ell_1^\Sigma(a/A)$ are domination-equivalent over $A$.
\item If $\tp(a/A)$ is $\Sigma\cup\Sigma'$-analysable and foreign to $\Sigma'$, then $a$ and $\ell_1^\Sigma(a/A)$ are domination-equivalent over $A$.\end{enumerate}
\esatz
\bew (1) Since $\ell_1^\Sigma(a/A)\in\bdd(Aa)$, clearly $a$ dominates (and $\Sigma$-dominates) $\ell_1^\Sigma(a/A)$ over $A$.

For the converse, suppose $\tp(b/A)$ is $\Sigma$-analysable and $b\nind_Aa$. Consider a sequence $(b_i:i<\alpha)$ in $\bdd(Ab)$ such that $\tp(b_i/A,b_j:j<i)$ is $\Sigma$-internal for all $i<\alpha$ and $b\in\bdd(A,b_i:i<\alpha)$. Since $a \nind_A b$ there is a minimal $i<\alpha$ such that $a\nind_{A,(b_j:j<i)}b_i$. Put 
$$a'=\cb(b_j:j\le i/Aa)\in\bdd(Aa).$$ 
Then $\tp(a'/A)$ is $\Sigma$-internal by Corollary \ref{cbase1}, and $a'\subseteq\ell_1^\Sigma(a/A)$.
Clearly $a'\nind_A (b_j:j\le i)$, whence $a'\nind_A b$ and finally $\ell_1^\Sigma(a/A)\nind_A b$. This shows $(1)$.

$(2)$ follows from $(3)$ setting $\Sigma'=\emptyset$.

$(3)$ Suppose $b\nind_A a$. We may assume that $b=\cb(a/Ab)$, so $\tp(b/A)$ is $(\Sigma\cup\Sigma')$-analysable. Thus 
$$b\nind_A \ell_1^{\Sigma\cup\Sigma'}(a/A)$$
by $(1)$. Now $\tp(\ell_1^{\Sigma\cup\Sigma'}(a/A)/A)$ is foreign to $\Sigma'$ since $\tp(a/A)$ is; it is hence almost $\Sigma$-internal. Therefore $\ell_1^{\Sigma\cup\Sigma'}(a/A)\subseteq\ell_1^\Sigma(a/A)$ and so $b\nind_A \ell_1^\Sigma(a/A)$.\qed

\bem If $\tp(a/A)$ is $\Sigma_0$-analysable and $\Sigma_1$ is a subfamily of $\Sigma_0$ such that $\tp(a/A)$ remains $\Sigma_1$-analysable, then
$$\ell_1^{\Sigma_1}(a/A)\subseteq\ell_1^{\Sigma_0}(a/A)\subseteq\bdd(aA)$$
and $\ell_1^{\Sigma_1}(a/A)$ and $\ell_1^{\Sigma_0}(a/A)$ are both domination-equivalent to $a$ over $A$.
In fact it would be sufficient to have $\Sigma_1$ such that $\tp(\ell_1^{\Sigma_0}(a/A)/A)$ is $\Sigma_1$-analysable.\ebem
\frag When is there a minimal (boundedly closed) $a_0\in\bdd(aA)$ domination-equivalent with $a$ over $A$?\efrag
If $T$ has finite $\SU$-rank, one can take $a_0\in\bdd(aA)\setminus\bdd(A)$ with $SU(a_0/A)$ minimal possible.

\defn\begin{itemize}\item A type $\tp(a/A)$ is {\em $\Sigma$-flat} if $\ell_1^\Sigma(a/A)=\ell_2^\Sigma(a/A)$. It is {\em $A$-flat} if it is $\Sigma$-flat for all $A$-invariant $\Sigma$. It is {\em flat} if for all $B\supseteq A$ every nonforking extension to $B$ is $B$-flat. A theory $T$ is {\em flat} if all its types are.
\item A type $p\in S(A)$ is {\em $A$-ultraflat} if it is almost internal in any $A$-invariant family of partial types it is non-foreign to. It is {\em ultraflat} if for any $B\supseteq A$ every nonforking extension to $B$ is $B$-ultraflat.\end{itemize}\edefn

Flatness and ultraflatness are clearly preserved under non-forking extensions and non-forking restrictions, and under adding and forgetting parameters.

\bem If $\tp(a/A)$ is $\Sigma$-flat, then $\ell_\alpha^\Sigma(a/A)=\ell_1^\Sigma(a/A)$ for all $\alpha>0$. Clearly, ultraflat implies flat.\ebem
\bsp\begin{itemize}
\item Generic types of fields or definably simple groups interpretable in a simple theory are ultraflat.
\item Types of Lascar rank $1$ are ultraflat.
\item If there is no boundedly closed set between $\bdd(A)$ and $\bdd(aA)$, then $\tp(a/A)$ is $A$-ultraflat.
\item In a small simple theory there are many $A$-ultraflat types over finite sets $A$, as the lattice of boundedly closed subsets of $\bdd(aA)$ is scattered for finitary $aA$.\end{itemize}\ebsp

Next we shall prove that any type internal in a family of Lascar rank one types is also flat.

\lmm\label{flat-bdd} It $\tp(a/A)$ is flat (ultraflat), then so is $\tp(a_0/A)$ for any $a_0\in\bdd(Aa)$.\elmm
\bew Consider a set $B$ extending $A$ with $B\ind_Aa_0$; we may assume that $B\ind_{Aa_0} a$,whence $B\ind_A a$. 

Firstly, the flat case is clear since $\ell_\alpha^\Sigma(a_0/B)=\ell_\alpha^\Sigma(a/B)\cap\bdd(Ba_0)$ for any $\alpha>0$ and for any $B$-invariant family $\Sigma$. Assume now $\tp(a/A)$ is ultraflat, $a_0\in\bdd(Aa)$ and $\tp(a_0/B)$ is not foreign to some $B$-invariant family $\Sigma$. Then $\tp(a/B)$ is not foreign to $\Sigma$, hence almost $\Sigma$-internal, as is $\tp(a_0/B)$. \qed

\kor If $\tp(a/A)$ is almost internal in a family of types of Lascar rank one, then it is flat.\ekor 
\bew Assume there is some $B\ind_A a$ and some tuple $\bar b$ of realizations of types of Lascar rank one over $B$ such that $a\subseteq\bdd(B\bar b)$. We may assume $\bar b$ is an independent sequence over $B$ since all its elements have $\SU$-rank one. Hence $\bar b$ is an independent sequence over any $C\supseteq B$ with $C\ind_B \bar b$, so $\tp(\bar b/B)$ is flat by Corollary \ref{cor-level1-ind}. Thus, $\tp(a/B)$ is flat by Lemma  \ref{flat-bdd}, and so is $\tp(a/A)$. \qed

\frag Is every (finitary) type in a small simple theory non-orthogonal to a flat type?\efrag
\frag Is every type in a supersimple theory non-orthogonal to a flat type?\efrag
\prob Construct a flat type which is not ultraflat.\eprob

We shall now recall the definitions and properties of $\Sigma$-closure from \cite[Section 4.0]{wa97} in the stable and \cite[Section 3.5]{wa00} in the simple case (where it is called $P$-closure: our $\Sigma$ corresponds to the collection of all $P$-analysable types which are co-foreign to $P$). Buechler and Hoover \cite[Definition 1.2]{bh01} redefine such a closure operator in the context of superstable theories and reprove some of the properties \cite[Lemma 2.5]{bh01}.
\defn The {\em $\Sigma$-closure\/} $\Pcl(A)$ of a set $A$ is the
collection of all hyperimaginaries $a$ such that $\tp(a/A)$ is
$\Sigma$-analysable.\edefn
\bem We think of partial types in $\Sigma$ as small. We always have $\bdd(A)\subseteq\Pcl(A)$; equality holds if
$\Sigma$ is the family of all bounded types. Other useful examples for $\Sigma$ are the family of all
types of $SU$-rank $<\omega^\alpha$ for some ordinal $\alpha$, the
family of all supersimple types in a properly simple theory, or the family of $p$-simple types of $p$-weight $0$ for some regular type $p$, giving rise to Hrushovski's $p$-closure \cite{udi1}.\ebem
\tats\label{forequ} The following are equivalent:\begin{enumerate}
\item $\tp(a/A)$ is foreign to $\Sigma$.
\item $a\ind_A\Pcl(A)$.
\item $a\ind_A\dcl(aA)\cap\Pcl(A)$.
\item $\dcl(aA)\cap\Pcl(A)\subseteq\bdd(A)$.\end{enumerate}\etats
\bew The equivalence of (1), (2) and (3) is \cite[Lemma 3.5.3]{wa00}; the equivalence $(3)\Leftrightarrow(4)$ is obvious.\qed

Unless it equals bounded closure, $\Sigma$-closure has the size of the monster model and thus violates the usual conventions. The equivalence $(2)\Leftrightarrow(3)$ can be used to cut it down to some small part.
\tats\label{clind} Suppose $A\ind_BC$. Then $\Pcl(A)\ind_{\Pcl(B)}\Pcl(C)$.
More precisely, for any $A_0\subseteq\Pcl(A)$ we have
$A_0\ind_{B_0}\Pcl(C)$, where
$B_0=\dcl(A_0B)\cap\Pcl(B)$. In particular,
$\Pcl(AB)\cap\Pcl(BC)=\Pcl(B)$.\etats
\bew This is \cite[Lemma 3.5.5]{wa00}; the second clause follows from Fact \ref{forequ}.\qed

\lmm\label{intersection} Suppose $C\subseteq A\cap B\cap D$ and $AB\ind_CD$.\begin{enumerate}
\item If $\Pcl(A)\cap\Pcl(B)=\Pcl(C)$, then $\Pcl(AD)\cap\Pcl(BD)=\Pcl(D)$.
\item If $\bdd(A)\cap\Pcl(B)=\bdd(C)$, then $\bdd(AD)\cap\Pcl(BD)=\bdd(D)$.
\end{enumerate}\elmm
\bew (1) This is \cite[Lemma 3.5.6]{wa00}, which in turn adapts \cite[Fact 2.4]{pi95}.

(2) This is similar to $(1)$. By Fact 15 $$\Pcl(BD)\ind_{\Pcl(B)\cap\dcl(AB)} AB\ ;$$
since $AD\ind_A AB$ we obtain
$$\cb(\bdd(AD)\cap\Pcl(BD)/AB)\subseteq \bdd(A)\cap\Pcl(B)=\bdd(C).$$ Hence
$$\bdd(AD)\cap\Pcl(BD)\ind_C AB$$ and by transitivity
$$\bdd(AD)\cap\Pcl(BD)\ind_D ABD,$$
whence the result.\qed

The following lemma tells us that we can actually find a set $C$ with $\Pcl(A)\cap\Pcl(B)=\Pcl(C)$ as in Lemma \ref{intersection}(1), even though the $\Sigma$-closures have the size of the monster model.
\lmm\label{Pintersection} Let $C=\bdd(AB)\cap\Pcl(A)\cap\Pcl(B)$. Then $\Pcl(A)\cap\Pcl(B)=\Pcl(C)$.\elmm
\bew Consider $e\in\Pcl(A)\cap\Pcl(B)$ and put $f=\cb(e/AB)$. Then $e\ind_fAB$~; since $\tp(e/A)$ is $\Sigma$-analysable, so is $\tp(e/f)$, and $e\in\Pcl(f)$. If $I$ is a Morley sequence in $\tp(e/AB)$, then $f\in\dcl(I)$. However, since $e$ is $\Sigma$-analysable over $A$ and over $B$, so is $I$, whence $f$. Hence
$$f\in\bdd(AB)\cap\Pcl(A)\cap\Pcl(B)=C.$$
The result follows.\qed

However, for considerations such as the canonical base property, one should like to work with the first level of the $\Sigma$-closure rather than with the full closure operator.
\defn The $\Sigma_1$-closure of $A$ is given by
$$\Qcl(A)=\ell_1^\Sigma(\Pcl(A)/A)=\{b:\tp(b/A)\mbox{ is almost $\Sigma$-internal}\}.$$\edefn
Unfortunately, unless $\tp(\Pcl(A)/A)$ is $\Sigma$-flat, $\Sigma_1$-closure is not a closure operator, as $\Qcl(\Qcl(A))\supset\Qcl(A)$.
\lmm\label{1clind} Suppose $A\ind_BC$ with $B\subseteq A\cap C$. Then
$$\Qcl(A)\ind_{\Qcl(B)}C.$$
More precisely, $\Qcl(A)\ind_{\Qcl(B)\cap\bdd(C)}C$.\elmm
\bew Consider $a\in\Qcl(A)$ and put $c=\cb(Aa/C)$. Then $\tp(c/B)$ is almost $\Sigma$-internal by Corollary \ref{cbase1}, and $c\in\bdd(C)\cap\Qcl(B)$.\qed
\frag If $A\ind_BC$, is $\Qcl(A)\ind_{\Qcl(B)}\Qcl(C)$~?\efrag

\section{$\Sigma$-ampleness and weak $\Sigma$-ampleness}
Let $\Phi$ and $\Sigma$ be $\emptyset$-invariant families of partial types.
\defn\label{ample} $\Phi$ is {\em $n$-$\Sigma$-ample} if there are tuples $a_0,\ldots,a_n$, with $a_n$ a tuple of realizations of partial types in $\Phi$ over some parameters $A$, such that\begin{enumerate}
\item $a_n\nind_{\Pcl(A)} a_0$;
\item $a_{i+1}\ind_{\Pcl(Aa_i)}a_0\ldots a_{i-1}$ for $1\le i<n$;
\item $\Pcl(Aa_0\ldots a_{i-1}a_i)\cap\Pcl(Aa_0\ldots a_{i-1}a_{i+1})=\Pcl(Aa_0\ldots a_{i-1})$ for $0\le i<n$.\end{enumerate}\edefn
\bem Pillay \cite{pi00} requires $a_n\ind_{Aa_i}a_0\ldots a_{i-1}$ for $1\le i<n$ in item (2). We follow the variant proposed by Evans and N\"ubling \cite{ev03} which seems more natural and which implies $$a_n\ldots a_{i+1}\ind_{\Pcl(Aa_i)}a_0\ldots a_{i-1}.$$\ebem
\lmm\label{monotone} If $\Sigma'$ is a $\Sigma$-analysable family of partial types, then $n$-$\Sigma$-ample implies $n$-$\Sigma'$-ample, and in particular $n$-ample.\elmm
\bew As in \cite[Remark 3.7]{pi00} we replace $a_i$ by 
$$a'_i=\cb(a'_n\ldots a'_{i+1}/\Pcl(Aa_i))$$
for $i<n$, where $a'_n=a_n$. Then
$$a'_n\ldots a'_{i+1}\ind_{a'_i}\Pcl(Aa_i)\quad\text{and}\quad a'_n\ldots a'_{i+1}\ind_{\Pcl(Aa_i)}\Pcl(Aa_0\ldots a_i)$$
by Fact \ref{forequ}, whence
$$a'_n\ldots a'_{i+1}\ind_{a'_i}\Pcl(Aa_0\ldots a_i).$$
Put $A'=\Pcl(A)\cap\bdd(Aa'_0)$. Then
$A\subseteq A'\subseteq\Pcl(A)$, whence $\Pcl(A)=\Pcl(A')$, and $a'_0\ind_{A'}\Pcl(A)$.
Now $a_n\nind_{\Pcl(A')}a_0$ implies $a'_n\nind_{\Pcl(A)}a'_0$, whence $a'_n\nind_{\PPcl(A)}a'_0$. Clearly $a'_{i+1}\ind_{a'_i}\Pcl(Aa_0\ldots a_i)$ implies $$a'_{i+1}\ind_{\PPcl(A'a'_i)}a'_0\ldots a'_{i-1}$$ for $1\le i<n$. Finally,
$$A'a'_0\ldots a'_ia'_{i+1}\ind_{\PPcl(A'a'_0\ldots a'_{i-1})}\Pcl(Aa_0\ldots a_{i-1})$$
yields
$$\PPcl(A'a'_0\ldots a'_{i-1}a'_{i+1}) ,\PPcl(A'a'_0\ldots a'_i)\ind_{\PPcl(A'a'_0\ldots a'_{i-1})}\Pcl(Aa_0\ldots a_{i-1}),$$
so
$$\PPcl(A'a'_0\ldots a'_{i-1}a'_{i+1})\cap\PPcl(A'a'_0\ldots a'_i)\subseteq\Pcl(Aa_0\ldots a_{i-1})$$
implies
$$\PPcl(A'a'_0\ldots a'_{i-1}a'_{i+1})\cap\PPcl(A'a'_0\ldots a'_i)\subseteq\PPcl(A'a'_0\ldots a'_{i-1}).$$\qed

This also shows that in Definition \ref{ample} one may require $a_0,\ldots,a_{n-1}$ to lie in $\Phi^{heq}$, and $a_{i+1}\ind_{a_i}\Pcl(Aa_0\ldots a_i)$.
\bem\cite[Lemma 3.2 and Corollary 3.3]{pi00} If $a_0,\ldots,a_n$ witness $n$-$\Sigma$-ampleness over $A$, then $a_n\nind_{\Pcl(Aa_0\ldots a_{i-1})}a_i$ for all $i<n$. Hence $a_i,\ldots,a_n$ witness $(n-i)$-$\Sigma$-ampleness over $Aa_0\ldots a_{i-1}$. Thus $n$-$\Sigma$-ample implies $i$-$\Sigma$-ample for all $i\le n$.\ebem
\bem It is clear from the definition that even though $\Phi$ might be a complete type $p$, if $p$ is not $n$-$\Sigma$-ample, neither is any extension of $p$, not only the non-forking ones.\ebem
For $n=1$ and $n=2$ there are alternative definitions of non-$n$-$\Sigma$-ampleness:
\defn\begin{enumerate}
\item $\Phi$ is {\em $\Sigma$-based} if $\cb(a/\Pcl(B))\subseteq\Pcl(aA)$ for any tuple $a$ of realizations of partial types in $\Phi$ over some parameters $A$ and any $B\supseteq A$.
\item $\Phi$ is {\em $\Sigma$-CM-trivial} if $\cb(a/\Pcl(AB))\subseteq\Pcl(A,\cb(a/\Pcl(AC))$
for any tuple $a$ of realizations of partial types in $\Phi$ over some parameters $A$ and any $B\subseteq C$ such that $\Pcl(ABa)\cap\Pcl(AC)=\Pcl(AB)$.\end{enumerate}\edefn
\lmm\label{equivalence}\begin{enumerate}
\item $\Phi$ is $\Sigma$-based if and only if $\Phi$ is not $1$-$\Sigma$-ample.
\item $\Phi$ is $\Sigma$-CM-trivial if and only if $\Phi$ is not $2$-$\Sigma$-ample.
\end{enumerate}\elmm
\bew (1) Suppose $\Phi$ is $\Sigma$-based and consider $a_0,a_1,A$ with $\Pcl(Aa_0)\cap\Pcl(Aa_1)=\Pcl(A)$. Put
$a=a_1$ and $B=Aa_0$. By $\Sigma$-basedness
$$\cb(a/\Pcl(B))\subseteq\Pcl(Aa)\cap\Pcl(B)=\Pcl(A).$$
Hence $a\ind_{\Pcl(A)}\Pcl(B)$, whence $a_1\ind_{\Pcl(A)}a_0$, so $\Phi$ is not $1$-$\Sigma$-ample.

Conversely, if $\Phi$ is not $\Sigma$-based, let $a,A,B$ be a counterexample. Put $a_0=\cb(a_1/\Pcl(B))$ and $a_1=a$. Then $a_0\notin\Pcl(Aa_1)$. Now take
$$A'=\bdd(Aa_0a_1)\cap\Pcl(Aa_0)\cap\Pcl(Aa_1).$$
Then $\Pcl(A'a_0)\cap\Pcl(A'a_1)=\Pcl(A')$ by Lemma \ref{Pintersection}.

Suppose $a_1\ind_{\Pcl(A')}a_0$. Since $\Pcl(A')\subseteq\Pcl(Aa_0)\subseteq\Pcl(B)$ we have $a_1\ind_{a_0}\Pcl(A')$. As $a_0=\cb(a_1/\Pcl(B))$, this implies
$$a_0\subseteq\Pcl(A')\subseteq\Pcl(Aa_1),$$
a contradiction. Hence $a_0,a_1,A'$ witness $1$-$\Sigma$-ampleness of $\Phi$.

(2) Suppose $\Phi$ is $\Sigma$-CM-trivial and consider $a_0,a_1,a_2,A$ with
$$\begin{aligned}\Pcl(Aa_0)\cap\Pcl(Aa_1)&=\Pcl(A),\\
\Pcl(Aa_0a_1)\cap\Pcl(Aa_0a_2)&=\Pcl(Aa_0),\quad\mbox{and}\\
a_2\ind_{\Pcl(Aa_1)}&a_0.\end{aligned}$$
Put $a=a_2$, $B=a_0$ and $C=a_0a_1$. Then
$$a_2\ind_{\Pcl(Aa_1)}\Pcl(Aa_0a_1),$$
so $\cb(a/\Pcl(AC))\subseteq\Pcl(Aa_1)$. Moreover
$$\Pcl(ABa)\cap\Pcl(AC)=\Pcl(AB),$$
whence by $\Sigma$-CM-triviality
$$\begin{aligned}\cb(a/\Pcl(AB))&\subseteq\Pcl(A,\cb(a/AC))\cap\Pcl(AB)\\
&\subseteq\Pcl(Aa_1)\cap\Pcl(Aa_0)=\Pcl(A).\end{aligned}$$
Hence $a_2\ind_{\Pcl(A)}a_0$, so $\Phi$ is not $2$-$\Sigma$-ample.

Conversely, if $\Phi$ is not $\Sigma$-CM-trivial, let $a,A,B,C$ be a counterexample. Put
$$a_0=\cb(a/\Pcl(AB)),\quad a_1=\cb(a/\Pcl(AC)),\quad a_2=a,$$
$$A'=\bdd(Aa_0a_1)\cap\Pcl(Aa_0)\cap\Pcl(Aa_1)\subseteq\Pcl(AB).$$
Then $a_2\ind_{\Pcl(A'a_1)}a_0$ and $a_0\notin\Pcl(Aa_1)$; by Lemma \ref{Pintersection}
$$\Pcl(A'a_0)\cap\Pcl(A'a_1)=\Pcl(A').$$
Moreover, $a_2\ind_{a_0}\Pcl(AB)$ implies
$$\Pcl(A'a_0a_2)\ind_{\Pcl(A'a_0)}\Pcl(AB).$$
Thus
$$\begin{aligned}\Pcl(A'a_0a_2)&\cap\Pcl(A'a_0a_1)\subseteq\Pcl(ABa)\cap\Pcl(AC)\\
&=\Pcl(AB)\cap\Pcl(A'a_0a_2)=\Pcl(A'a_0).\end{aligned}$$
Suppose $a_2\ind_{\Pcl(A')}a_0$. Since $a_2\ind_{a_0}\Pcl(A')$ we obtain
$$a_0=\cb(a/\Pcl(AB))=\cb(a/a_0\Pcl(A'))\subseteq\Pcl(A')\subseteq\Pcl(Aa_1),$$
a contradiction. Hence $a_0,a_1,a_2,A'$ witness $2$-$\Sigma$-ampleness of~$\Phi$.\qed

In our definition of $\Sigma$-ampleness, we only consider the type of $a_n$ over a $\Sigma$-closed set, namely $\Pcl(A)$. This seems natural since the idea of $\Sigma$-closure is to work {\em modulo $\Sigma$}. However, sometimes one needs a stronger notion which takes care of all types. Let us first look at $n=1$ and $n=2$.
\defn\begin{itemize}\item $\Phi$ is {\em strongly $\Sigma$-based} if $\cb(a/B)\subseteq\Pcl(aA)$ for any tuple $a$ of realizations of partial types in $\Phi$ over some $A$ and any $B\supseteq A$.
\item $\Phi$ is {\em strongly $\Sigma$-CM-trivial} if $\cb(a/AB)\subseteq\Pcl(A,\cb(a/AC)$ for any tuple $a$ of realizations of partial types in $\Phi$ over some $A$ and any $B\subseteq C$ with $\Pcl(ABa)\cap\bdd(AC)=\bdd(AB)$.\end{itemize}\edefn

\bem $\cb(a/\Pcl(B))\subseteq\bdd(\cb(a/B),a)\cap\Pcl(\cb(a/B))$. \ebem
\bew By Fact \ref{clind} the independence $a\ind_{\cb(a/B)}B$ implies
$$a\ind_{\dcl(a,\cb(a/b))\cap\Pcl(\cb(a/B))}\Pcl(B).$$
The result follows.\qed

\verm $\cb(a/B)\subseteq\Pcl(\cb(a/\Pcl(B)))$.\everm
If this conjecture were true, strong and ordinary $\Sigma$-basedness and $\Sigma$-CM-triviality would obviously coincide. Since we have not been able to show this, we weaken our definition of ampleness.

\defn\label{wample} $\Phi$ is {\em weakly $n$-$\Sigma$-ample} if there are tuples $a_0,\ldots,a_n$, where $a_n$ is a tuple of realizations of partial types in $\Phi$ over $A$, with\begin{enumerate}
\item $a_n\nind_A a_0$.
\item $a_{i+1}\ind_{Aa_i}a_0\ldots a_{i-1}$ for $1\le i<n$.
\item $\bdd(Aa_0\ldots a_{i-1}a_i)\cap\Pcl(Aa_0\ldots a_{i-1}a_{i+1})=\bdd(Aa_0\ldots a_{i-1})$ for $i<n$.\end{enumerate}\edefn
Note that (3) implies that $\tp(a_i/Aa_0\ldots a_{i-1})$ is foreign to $\Sigma$ by Fact \ref{forequ} for all $i<n$, and so is $\tp(a_i/Aa_{i-1})$ by (2). If $\Sigma$ is the family of bounded partial types, then weak and ordinary $n$-$\Sigma$-ampleness just equal $n$-ampleness.

\lmm\label{ample-weak-ordinary} An $n$-$\Sigma$-ample family of types is weakly $n$-$\Sigma$-ample. If $\Sigma'$ is $\Sigma$-analysable, then a weakly $n$-$\Sigma$-ample family is weakly $n$-$\Sigma'$-ample, and in particular $n$-ample.\elmm
\bew If $a_0,\ldots,a_n$ witness $n$-$\Sigma$-ampleness over $A$, we put $a'_n=a_n$, 
$$a'_i=\cb(a'_n\ldots a'_{i+1}/\Pcl(Aa_i))\subseteq\Pcl(Aa_i)\quad\text{for $n>i$}$$
and 
$$A'=\bdd(Aa'_0)\cap\Pcl(Aa'_1)\subseteq\Pcl(Aa_0)\cap\Pcl(Aa_1)=\Pcl(A).$$
As in Lemma \ref{monotone} we have for $i<n$
$$a'_n\ldots a'_{i+1}\ind_{a'_i}\Pcl(Aa_0\ldots a_i).$$
For $0<i<n$ we obtain $a'_{i+1}\ind_{A'a'_i}a'_0\ldots a'_{i-1}$; moreover
$$\begin{aligned}\bdd(A'a'_0\ldots a'_{i-1}a'_i)&\cap\Pcl(A'a'_0\ldots a'_{i-1}a'_{i+1})\\
&\subseteq\Pcl(A'a'_0\ldots a'_{i-1}a'_i)\cap\Pcl(A'a'_0\ldots a'_{i-1}a'_{i+1})\\
&\subseteq\Pcl(Aa_0\ldots a_{i-1}a_i)\cap\Pcl(Aa_0\ldots a_{i-1}a_{i+1})\\
&=\Pcl(Aa_0\ldots a_{i-1}).\end{aligned}$$
But then $a'_i\ind_{A'a'_0\ldots a'_{i-1}}\Pcl(Aa_0\ldots a_{i-1})$ yields
$$\bdd(A'a'_0\ldots a'_{i-1}a'_i)\cap\Pcl(A'a'_0\ldots a'_{i-1}a'_{i+1})=\bdd(A'a'_0\ldots a'_{i-1}),$$
while $\bdd(A'a'_0)\cap\Pcl(A'a'_1)=\bdd(A')$ follows from the definition of
$A'$. Finally $a_n\nind_{\Pcl(A)}a_0$ implies $a'_n\nind_{\Pcl(A)}a'_0$, 
whence $a'_n\nind_{A'}a'_0$ as $\tp(a'_0/A')$ is foreign to $\Sigma$ and
$\Pcl(A)=\Pcl(A')$.

The second assertion is clear, since $\PPcl(A)\subseteq\Pcl(A)$ for any set $A$.\qed

This also shows that in Definition \ref{wample} one may require $a_0,\ldots,a_{n-1}$ to lie in $\Phi^{heq}$.

\lmm\begin{enumerate}
\item $\Phi$ is strongly $\Sigma$-based iff $\Phi$ is not weakly $1$-$\Sigma$-ample.
\item $\Phi$ is strongly $\Sigma$-CM-trivial iff $\Phi$ is not weakly $2$-$\Sigma$-ample.
\end{enumerate}\elmm
\bew This is similar to the proof of Lemma \ref{equivalence}, so we shall be concise.

(1) Suppose $\Phi$ is strongly $\Sigma$-based and consider $a_0,a_1,A$ with $$\bdd(Aa_0)\cap\Pcl(Aa_1)=\bdd(A).$$
Put $a=a_1$ and $B=Aa_0$. By strong $\Sigma$-basedness
$$\cb(a/B)\subseteq\Pcl(Aa)\cap\bdd(B)=\bdd(A),$$
whence $a_1\ind_Aa_0$, so $\Phi$ is not weakly $1$-$\Sigma$-ample.

Conversely, if $\Phi$ is not strongly $\Sigma$-based, let $a,A,B$ be a counterexample. Put $a_0=\cb(a_1/B)$ and $a_1=a$. Then $a_0\notin\Pcl(Aa_1)$. Now take $A'=\bdd(Aa_0)\cap\Pcl(Aa_1)$. Clearly $A'=\bdd(A'a_0)\cap\Pcl(A'a_1)$.
Suppose $a_1\ind_{A'}a_0$. Since $a_0=\cb(a_1/B)$ implies $a_1\ind_{a_0}A'$, we obtain
$$a_0\subseteq\bdd(A')\subseteq\Pcl(Aa_1),$$
a contradiction. Hence $a_0,a_1,A'$ witness weak $1$-$\Sigma$-ampleness of~$\Phi$.

(2) Suppose $\Phi$ is strongly $\Sigma$-CM-trivial and consider $a_0,a_1,a_2,A$ with
$$\begin{aligned}\bdd(Aa_0)\cap\Pcl(Aa_1)&=\bdd(A),\\
\bdd(Aa_0a_1)\cap\Pcl(Aa_0a_2)&=\bdd(Aa_0),\quad\mbox{and}\\
a_2\ind_{Aa_1}&a_0.\end{aligned}$$
Put $a=a_2$, $B=a_0$ and $C=a_0a_1$. Then $\cb(a/AC)\subseteq\bdd(Aa_1)$. Moreover
$$\Pcl(ABa)\cap\bdd(AC)=\bdd(AB),$$
whence by strong $\Sigma$-CM-triviality
$$\begin{aligned}\cb(a/AB)&\subseteq\Pcl(A,\cb(a/AC))\cap\bdd(AB)\\
&\subseteq\Pcl(Aa_1)\cap\bdd(Aa_0)=\bdd(A).\end{aligned}$$
Hence $a_2\ind_Aa_0$, so $\Phi$ is not $2$-$\Sigma$-ample.

Conversely, if $\Phi$ is not strongly $\Sigma$-CM-trivial, let $a,A,B,C$ be a counterexample. Put
$$a_0=AB,\quad a_1=\cb(a/AC),\quad a_2=a,$$
$$A'=\bdd(Aa_0)\cap\Pcl(Aa_1).$$
Then $a_2\ind_{A'a_1}a_0$ and $\cb(a_2/AB)\notin\Pcl(Aa_1)=\Pcl(A'a_1)$; moreover
$$\bdd(A'a_0)\cap\Pcl(A'a_1)=\bdd(A').$$
Clearly
$$\begin{aligned}\Pcl(A'a_0a_2)&\cap\bdd(A'a_0a_1)\subseteq\Pcl(ABa)\cap\bdd(AC)\\
&=\bdd(AB)=\bdd(A'a_0).\end{aligned}$$
Suppose $a_2\ind_{A'}a_0$. Then $\cb(a_2/AB)\in\bdd(A')\subseteq\Pcl(Aa_1)$,
a contradiction. Hence $a_0,a_1,a_2,A'$ witness weak $2$-$\Sigma$-ampleness of~$\Phi$.\qed

\lmm\label{clbase} If $\Phi$ is not (weakly) $n$-$\Sigma$-ample, neither is the family of $\emptyset$-conjugates of 
$\tp(a/A)$ for any $a\in\Pcl(\bar aA)$, where $\bar a$ is a tuple of realizations of partial types in $\Phi$ over $A$.\elmm
\bew Suppose the family of $\emptyset$-conjugates of $\tp(a/A)$ is $n$-$\Sigma$-ample, as witnessed by $a_0,\ldots,a_n$ 
over some parameters $B$. There is a tuple $\bar a$ of realizations of partial types in $\Phi$ over some 
$\emptyset$-conjugates of $A$ inside $B$ such that $a_n\in\Pcl(\bar aB)$; we may choose it such that
$$\bar a\ind_{a_nB}a_0\ldots a_{n-1}.$$
Then $\bar a\ind_{a_{n-1}a_nB}a_0\ldots a_{n-2}$, and hence
$$\bar a\ind_{\Pcl(a_{n-1}a_nB)}a_0\ldots a_{n-2}.$$
As $a_n\ind_{\Pcl(a_{n-1}B)}a_0\ldots a_{n-2}$ implies
$$\Pcl(a_{n-1}a_nB)\ind_{\Pcl(a_{n-1}B)}a_0\ldots a_{n-2}$$
by Fact \ref{clind}, we get
$$\bar a\ind_{\Pcl(a_{n-1}B)}a_0\ldots a_{n-2}.$$
We also have $\bar a\ind_{a_0\ldots a_{n-2}a_nB}a_{n-1}$, whence
\begin{equation}\label{eqn1}\Pcl(a_0\ldots a_{n-2}\bar aB)\ind_{\Pcl(a_0\ldots a_{n-2}a_nB)}\Pcl(a_0\ldots a_{n-2}a_{n-1}B);\end{equation}
since $\Sigma$-closure is boundedly closed,
$$\begin{aligned}\Pcl(a_0\ldots a_{n-2}\bar aB)&\cap\Pcl(a_0\ldots a_{n-2}a_{n-1}B)\\
&\subseteq\Pcl(a_0\ldots a_{n-2}a_nB)\cap\Pcl(a_0\ldots a_{n-2}a_{n-1}B)\\
&=\Pcl(a_0\ldots a_{n-2}B).\end{aligned}$$
Finally, $\bar a\ind_{\Pcl(B)}a_0$ would imply $\Pcl(\bar aB)\ind_{\Pcl(B)}a_0$ by Fact \ref{clind}, and hence $a_n\ind_{\Pcl(B)}a_0$, a contradiction. Thus $\bar a\nind_{\Pcl(B)}a_0$, and $a_0,\ldots,a_{n-1},\bar a$ witness $n$-$\Sigma$-ampleness of $\Phi$ over $B$, a contradiction.

Now suppose $a_0,\ldots,a_n$ witness weak $n$-$\Sigma$-ampleness over $B$, and choose $\bar a$ as before. Then easily $\bar aa_n\ind_{Ba_{n-1}}a_0\ldots a_{n-2}$, yielding (2) from the definition. Moreover, equation (\ref{eqn1}) implies
$$\begin{aligned}\Pcl(a_0\ldots a_{n-2}\bar aB)&\cap\bdd(a_0\ldots a_{n-2}a_{n-1}B)\\
&\subseteq\Pcl(a_0\ldots a_{n-2}a_nB)\cap\bdd(a_0\ldots a_{n-2}a_{n-1}B)\\
&=\bdd(a_0\ldots a_{n-2}B).\end{aligned}$$
Finally suppose $\bar a\ind_B a_0$. Since $\tp(a_0/B)$ is foreign to $\Sigma$, so is $\tp(a_0/B\bar a)$. Then $a_0\ind_{B\bar a}\Pcl(B\bar a)$ by Fact \ref{forequ}, whence $a_0\ind_B a_n$, a contradiction. Thus $\bar a\nind_Ba_0$, and $a_0,\ldots,a_{n-1},\bar a$ witness weak $n$-$\Sigma$-ampleness of $\Phi$ over $B$, again a contradiction.\qed

\lmm\label{addparams} Suppose $B\ind_Aa_0\ldots a_n$. If $a_0,\ldots,a_n$ witness (weak) $n$-$\Sigma$-ampleness over $A$, they witness (weak) $n$-$\Sigma$-ampleness over $B$.\elmm
\bew Clearly $B\ind_{a_0\ldots a_{i-1}A}a_0\ldots a_{i+1}A$, so Lemma \ref{intersection} yields
$$\Pcl(Ba_0\ldots a_{i-1}a_i)\cap\Pcl(Ba_0\ldots a_{i-1}a_{i+1})=\Pcl(Ba_0\ldots a_{i-1})$$ in the ordinary case, and
$$\bdd(Ba_0\ldots a_{i-1}a_i)\cap\Pcl(Ba_0\ldots a_{i-1}a_{i+1})=\bdd(Ba_0\ldots a_{i-1})$$ 
in the weak case, for all $i<n$.

Next, $a_{i+1}\ind_{Aa_0\ldots a_i}B$, whence $a_{i+1}\ind_{\Pcl(Aa_0\ldots a_i)}\Pcl(Ba_i)$ by Lemma \ref{clind}. Now $a_{i+1}\ind_{\Pcl(Aa_i)}a_0\ldots a_{i-1}$ implies $a_{i+1}\ind_{\Pcl(Aa_i)}\Pcl(Aa_0\ldots a_i)$, whence
$$a_{i+1}\ind_{\Pcl(B a_i)}a_0\ldots a_{i-1}$$
for $1\le i<n$ by transitivity. In the weak case, $a_{i+1}\ind_{Aa_i}a_0\ldots a_{i-1}$ implies $a_{i+1}\ind_{Aa_i}Ba_0\ldots a_{i-1}$ by transitivity, whence $a_{i+1}\ind_{Ba_i}a_0\ldots a_{i-1}$.

Finally, $a_n\ind_{\Pcl(A)}\Pcl(B)$ by Fact \ref{clind}, so $a_n\ind_{\Pcl(B)} a_0$ would imply $a_n\ind_{\Pcl(A)}a_0$, a contradiction. Hence $a_n\nind_{\Pcl(B)} a_0$. In the weak case, $a_n\ind_AB$ and $a_n\nind_Aa_0$ yield directly $a_n\nind_Ba_0$.\qed

\lmm\label{unionbase} Let $\Psi$ be an $\emptyset$-invariant family of types. If $\Phi$ and $\Psi$ are not (weakly) $n$-$\Sigma$-ample, neither is $\Phi\cup\Psi$.\elmm
\bew Suppose $\Phi\cup\Psi$ is weakly $n$-$\Sigma$-ample, as witnessed by $a_0,\ldots,a_n=bc$ over some parameters $A$, where $b$ and $c$ are tuples of realizations of partial types in $\Phi$
and $\Psi$, respectively. As $\Psi$ is not $n$-$\Sigma$-ample, we must have $c\ind_Aa_0$. Put $a_0'=\cb(bc/a_0A)$. Then  $\tp(a'_0/A)$ is internal in $\tp(b/A)$ by Corollary \ref{cbase1}. Put
$$a_n'=\cb(a'_0/a_nA).$$
Then $\tp(a'_n/A)$ is $\tp(a'_0/A)$-internal and hence $\tp(b/A)$-internal. Note that $a_n\nind_Aa_0$ implies $a_n\nind_Aa'_0$, whence
$$a'_n\nind_Aa'_0\quad\mbox{and}\quad a'_n\nind_Aa_0.$$
Moreover $a'_n\in\bdd(Aa_n)$, so $a_0,\ldots,a_{n-1},a'_n$ witness weak $n$-$\Sigma$-ample\-ness over $A$.

As $\tp(a'_n/A)$ is $\tp(b/A)$-internal, there is $B\ind_Aa'_n$ and a tuple $\bar b$ of realizations of $\tp(b/A)$ with $a'_n\in\dcl(B\bar b)$. We may assume
$$B\ind_{Aa'_n}a_0\ldots a_{n-1},$$
whence $B\ind_Aa_0\ldots a_{n-1}a'_n$. Hence $a_0,\ldots,a_{n-1},a'_n$ witness weak $n$-$\Sigma$-ampleness over $B$ by Lemma \ref{addparams}. As $a'_n\in\dcl(B\bar b)$, this contradicts non weak $n$-$\Sigma$-ampleness of $\Phi$ by Lemma \ref{clbase}.

The proof in the ordinary case is analogous, replacing $A$ by $\Pcl(A)$.\qed

\kor\label{limit} For $i<\alpha$ let $\Phi_i$ be an $\emptyset$-invariant family of partial types. If $\Phi_i$ is not (weakly) $n$-$\Sigma$-ample for all $i<\alpha$, neither is
$\bigcup_{i<\alpha}\Phi_i$.\ekor
\bew This just follows from the local character of forking and Lemma \ref{unionbase}.\qed

\lmm\label{indbase} If the family of $\emptyset$-conjugates of $\tp(a/A)$ is not (weakly) $n$-$\Sigma$-ample and $a\ind A$, then $\tp(a)$ is not (weakly) $n$-$\Sigma$-ample.\elmm
\bew Suppose $\tp(a)$ is (weakly) $n$-$\Sigma$-ample, as witnessed by $a_0,\ldots,a_n$ over some parameters $B$, where $a_n=(b_i:i<k)$ is a tuple of realizations of $\tp(a)$. For each $i<k$ choose $B_i\ind_{b_i} (B,a_0\ldots a_n,B_j:j<i)$ with $B_ib_i\equiv Aa$. Then $B_i\ind b_i$, whence $(B_i:i<k)\ind Ba_0\ldots a_n$. Then $a_0,\ldots,a_n$ witness (weak) $n$-$\Sigma$-ampleness over $(B,B_i:i<k)$ by Lemma \ref{addparams}, a contradiction, since $\tp(b_i/B_i)$ is an $\emptyset$-conjugate of $\tp(a/A)$ for all $i<k$.\qed

\bem In fact, in the above Lemma it suffices to merely assume that the single type $\tp(a/A)$ is not (weakly) $n$-$\Sigma$-ample in the theory $T(A)$, using Corollary \ref{limit}. It follows that ampleness is preserved under adding and forgetting parameters.\ebem

\kor\label{internal} Let $\Psi$ be an $\emptyset$-invariant family of types. If $\Psi$ is $\Phi$-internal and $\Phi$ is not (weakly) $n$-$\Sigma$-ample, neither is $\Psi$.\ekor
\bew Immediate from Lemmas \ref{clbase} and \ref{indbase}.\qed

\satz\label{analysable} Let $\Psi$ be an $\emptyset$-invariant family of types. If $\Psi$ is $\Phi$-analysable and $\Phi$ is not (weakly) $n$-$\Sigma$-ample, neither is $\Psi$.\esatz
\bew Suppose $\Psi$ is $n$-$\Sigma$-ample, as witnessed by $a_0,\ldots,a_n$ over some parameters $A$, where $a_n$ is a tuple of realizations of $\Psi$. Put $a'_n=\ell_1^\Phi(a_n/\Pcl(A)\cap\bdd(Aa_n))$. Then $a_n$ and $a'_n$ are domination-equivalent over $\Pcl(A)\cap\bdd(Aa_n)$ by Theorem \ref{domequ}; moreover $a_n$ and hence $a'_n$ are independent of $\Pcl(A)$ over $\Pcl(A)\cap\bdd(Aa_n)$ by Fact \ref{forequ}, so $a_n$ and $a'_n$ are domination-equivalent over $\Pcl(A)$. Thus $a_0,\ldots,a'_n$ witness non-$\Sigma$-ampleness over $A$, contradicting Corollary \ref{internal}.

For the weak case we put $a'_n=\ell_1^\Phi(a_n/A)$. So $a_n$ and $a'_n$ are domination-equivalent over $A$, whence $a'_n\nind_Aa_0$. Thus $a_0,\ldots,a'_n$ witness weak non-$\Sigma$-ampleness over $A$, contradicting again Corollary \ref{internal}.\qed

\section{Analysability of canonical bases}

As an immediate Corollary to Theorem \ref{analysable}, we obtain the following:
\satz\label{ampleana} Suppose every type in $T$ is non-orthogonal to a regular type, and let $\Sigma$ be the family of all $n$-ample regular types. Then $T$ is not weakly $n$-$\Sigma$-ample.\esatz
\bew A non $n$-ample type is not weakly $\Sigma$-ample by Lemma \ref{ample-weak-ordinary}. So all regular types are not weakly $n$-$\Sigma$-ample. But every type is analysable in regular types by the non-orthogonality hypothesis.\qed
\kor\label{sigmabased} Suppose every type in $T$ is non-orthogonal to a regular type. Then $\tp(\cb(a/b)/a)$ is analysable in the family of all non one-based regular types, for all $a$, $b$.\ekor
\bew This is just Theorem \ref{ampleana} for $n=1$.\qed

Note that a forking extension of a non one-based regular type of infinite rank may be one-based.

{\sloppy\bem In fact, the proof shows more. For every type $p$ let $\Sigma(p)$ be the collection of types in $\Sigma$ not foreign to $p$. Then $\tp(\cb(a/b)/a)$ is analysable in $\Sigma(\tp(\cb(a/b)))$. In particular, if $\tp(a)$ or $\tp(b)$ has rank less than $\omega^\alpha$, so does $\tp(\cb(a/b))$. Hence $\tp(\cb(a/b)/a)$ is analysable in the family of all non one-based regular types of rank less than $\omega^\alpha$.\ebem}

Corollary \ref{sigmabased} is due to Zo\'e Chatzidakis for types of finite $\SU$-rank in simple theories \cite[Proposition 1.10]{zoe}. In fact, she even obtains $\tp(\cb(a/b)/\bdd(a)\cap\bdd(b))$ to be analysable in the family of non one-based types of rank $1$, and to decompose in components, each of which is analysable in a non-orthogonality class of regular types. In infinite rank, one has to work modulo types of smaller rank. So let $\Sigma_\alpha$ be the collection of all partial types of $\SU$-rank $<\omega^\alpha$, and $\P$ be the family of non $\Sigma_\alpha$-based types of $\SU$-rank $\omega^\alpha$. Note that the $\Sigma_\alpha$-based types of $\SU$-rank $\omega^\alpha$ are precisely the locally modular types of $\SU$-rank $\omega^\alpha$.
\satz Let $b=\cb(a/\Acl(b))$ be such that $\SU(b)<\omega^{\alpha+1}$ for some ordinal $\alpha$ and let $A=\Acl(b)\cap\Acl(a)$. Then $\tp(b/A)$ is $(\Sigma_\alpha\cup\P)$-analysable.
\esatz
\bew Firstly, if $a\in\Acl(b)$ then $a=b\in A$. Similarly, if $b\in\Acl(a)$ then $b\in A$; in both cases $\tp(b/A)$ is trivially $(\Sigma_\alpha\cup\P)$-analysable. Hence we may assume $a\not\in\Acl(b)$ and $b\not\in\Acl(a)$.

Suppose towards a contradiction that the result is false and consider a counterexample $a,b$ with $\SU(b)$ minimal modulo $\omega^\alpha$ and then $\SU(b/\Acl(a))$ being maximal modulo $\omega^\alpha$. Note that this implies $$\omega^\alpha\le\SU(b/a)\le\SU(b/A)\le\SU(b)<\omega^{\alpha+1}.$$
Clearly (after adding parameters) we may assume $A=\Acl(\emptyset)$. Then for any $c$ the type $\tp(c)$ is $(\Sigma_\alpha\cup\P$)-analysable iff $\tp(c/A)$ is.
\beh We may assume $a=\cb(b/\Acl(a))$.\ebeh
\bewbeh Put $\tilde a=\cb(b/\Acl(a))$ and $\tilde b=\cb(\tilde a/\Acl(b))$. Then $\tilde a\in\Acl(a)$ and $a\ind_{\tilde a}b$. Hence $\Acl(b)=\Acl(\tilde b)$ by \cite[Lemma 3.5.8]{wa00}, and $\tp(\tilde b)$ is not $(\Sigma_\alpha\cup\P)$-analysable either. Thus the pair $\tilde a,\tilde b$ also forms a counterexample. Moreover, $\SU(b)$ equals $\SU(\tilde b)$ modulo $\omega^\alpha$ and $\SU(b/\Acl(a))=\SU(b/\Acl(\tilde a))$ equals $\SU(\tilde b/\Acl(\tilde a))$ modulo $\omega^\alpha$.\qed

Since $a$ is definable over a finite part of a Morley sequence in $\lstp(b/a)$ by supersimplicity of $\tp(b)$, we see that $\SU(a)<\omega^{\alpha+1}$. On the other hand, $a\notin\Acl(b)$ implies $\SU(a/b)\ge\omega^\alpha$.

Let $\hat a\subseteq\bdd(a)$ and $\hat b\subseteq\bdd(b)$ be maximal $(\Sigma_\alpha\cup\P)$-analysable. Then $a\notin\Acl(\hat a)$ and $b\not\in\Acl(\hat b)$, and $\tp(a/\hat a)$ and $\tp(b/\hat b)$ are foreign to $\Sigma_\alpha\cup\P$. Since $\cb(\hat a/b)$ and $\cb(\hat b/a)$ are $(\Sigma_\alpha\cup\P)$-analysable, we obtain
$$a\ind_{\hat a}\hat b\quad\mbox{and}\quad b\ind_{\hat b}\hat a.$$
\beh\label{claim1} $\tp(b/\hat b)$ and $\tp(a/\hat a)$ are both $\Sigma_\alpha$-based.\ebeh
\bewbeh Let $\Phi$ be the family of $\Sigma_\alpha$-based types of $\SU$-rank $\omega^\alpha$. Then $\tp(a/\hat a)$ is $(\Sigma_\alpha\cup\P\cup\Phi)$-analysable, but foreign to $\Sigma_\alpha\cup\P$. Put $a_0=\ell_1^\Phi(a/\hat a)$ and $b_0=\ell_1^\Phi(b/\hat b)$. Then $a\domeq_{\hat a}a_0$ and $b\domeq_{\hat b}b_0$ by Lemma \ref{domequ}(3); as $a\ind_{\hat a}\hat b$ and $b\ind_{\hat b}\hat a$ we even have $a\domeq_{\hat a\hat b}a_0$ and $b\domeq_{\hat a\hat b}b_0$. Since $a\nind_{\hat a\hat b}b$ we obtain $a_0\nind_{\hat a\hat b}b_0$. Moreover, $\tp(a_0/\hat a)$ and $\tp(b_0/\hat b)$ are $\Sigma_\alpha$-based by Theorem \ref{analysable} (or \cite[Theorem 11]{wa04}).

On the other hand, as $a_0\nind_{\hat b}b_0$, we see that $b'=\cb(a_0/\Acl(b_0))$ is not contained in $\hat b$ and hence is not $(\Sigma_\alpha\cup\P)$-analysable. So $a_0,b'$ is another counterexample; by minimality of $\SU$-rank $b$ and $b'$ have the same $\SU$-rank modulo $\omega^\alpha$, whence $b\in\Acl(b_0)$. Hence $\tp(b/\hat b)$ is $\Sigma_\alpha$-based, as is $\tp(a/\hat a)$ since $a=\cb(b/a)$ and $a\ind_{\hat a}\hat b$.\qed

\beh\label{claim2} $\Acl(a,\hat b)=\Acl(b,\hat a)=\Acl(a,b)$.\ebeh
\bewbeh As $\tp(a/\hat a)$ is $\Sigma_\alpha$-based, we have
$$a\ind_{\Acl(a)\cap\Acl(\hat a b)} \hat a b,$$
whence
$$\Acl(a)\ind_{\Acl(a)\cap\Acl(\hat a b)} b$$
by Fact \ref{clind}. Thus $a=\cb(b/\Acl(a))\in\Acl(\hat a b)$. Similarly $b\in\Acl(\hat b a)$.\qed

Let now $(b)^\frown(b_j:j<\omega)$ be a Morley sequence in $\tp(b/a)$ and let $\hat b_j$ represent the part of $b_j$ corresponding to $\hat b$. Then $(\hat b_j:j<\omega)$ is a Morley sequence in $\tp(\hat b/\hat a)$ since $a\ind_{\hat a} \hat b$. As $\SU(\hat b)<\infty$ there is some minimal $m<\omega$ such that $\hat a=\cb(\hat b/\hat a)\in\Acl(\hat b,\hat b_j:j< m)$. Then $m>0$, as otherwise $\Acl(b)=\Acl(\hat a,b)\ni a$, which is impossible. Moreover, $a\in\Acl(\hat a,b_j)$ for all $j<m$ by invariance and hence, $a\in\Acl(\hat b,b_j:j<m)$.

Put $b'=\cb(b_j:j<m/\Acl(b))$. Then $(b_j:j<m)\ind_{b'\hat b} \Acl(b)$, so by Fact \ref{clind}
$$\Acl(\hat b,b_j:j<m)\ind_{\Acl(b',\hat b)} \Acl(b).$$
Then $a\ind_{\Acl(b',\hat b)} \Acl(b)$, so $b=\cb(a/\Acl(b))\in\Acl(b',\hat b)$. As $b\not\in\Acl(\hat b)$ we obtain $b'\not\in\Acl(\hat b)$.\
\beh $\tp(b'/\Acl(b')\cap\Acl(b_j:j<m))$ is not $(\Sigma_\alpha\cup\P)$-analysable.\ebeh
\bewbeh Note first that $(b_j:j<m)\ind_ab$ implies
$$\Acl(b_j:j<m)\ind_{\Acl(a)} \Acl(b)$$ by Fact \ref{clind}, whence
$$\Acl(b')\cap\Acl(b_j:j<m)\subseteq\Acl(b)\cap\Acl(a)=\Acl(\emptyset).$$
As $b\in\Acl(b',\hat b)$ and $\tp(b/\hat b)$ is not $(\Sigma_\alpha\cup\P)$-analysable, neither is $\tp(b'/\hat b)$, nor {\em a fortiori} $\tp(b'/\Acl(\emptyset))$.\qed

As $b'=\cb(b_j:j<m/\Acl(b'))$, the pair $(b_j:j<m),b'$ forms another counterexample. By minimality $\SU(b)$ equals $\SU(b')$ modulo $\omega^\alpha$, which implies $\Acl(b)=\Acl(b')$.

As $\tp(b_j/\hat b_j)$ is foreign to $\Sigma_\alpha\cup\P$ and $\hat b$ is $(\Sigma_\alpha\cup\P)$-analysable, we obtain
$\hat b\ind_{(\hat b_j:j<m)}(b_j:j<m)$ and hence by Fact \ref{clind}
$$\hat b\ind_{\Acl(\hat b_j:j<m)} \Acl(b_j:j<m).$$
On the other hand, as $\hat a\in\Acl(\hat b,\hat b_j:j<m)$ but $\hat a\not\in\Acl(\hat b_j:j<m)$ by minimality of $m$, we get
$$\SU(\hat b/\Acl(\hat b_j:j<m)>_\alpha \SU(\hat b/\hat a,\Acl(\hat b_j:j<m)),$$
where the index $\alpha$ indicates modulo $\omega^\alpha$.

Moreover, as $\hat b\ind_{\hat a} a$ we get $\hat b\ind_{\Acl(\hat a)} \Acl(a)$, i.e.\ $\SU(\hat b/\Acl(\hat a))=\SU(\hat b/\Acl(a))$. Since $\Acl(b)=\Acl(b')$ and $b\in\Acl(a\hat b)$ we obtain
\begin{align*} \SU(b'/&\Acl(b_j:j<m))=_\alpha \SU(b/\Acl(b_j:j<m))\\
&\geq_\alpha \SU(\hat b/\Acl(b_j:j<m))=_\alpha \SU(\hat b/\Acl(\hat b_j:j<m))\\
&>_\alpha \SU(\hat b/\hat a,\Acl(\hat b_j:j<m))=_\alpha \SU(\hat b/\Acl(\hat a))\\
&=_\alpha\SU(\hat b/\Acl(a))=_\alpha\SU(b/\Acl(a)),\end{align*}
contradicting the maximality of $\SU(b/\Acl(a))$ modulo $\omega^\alpha$. This finishes the proof.\qed

As a corollary we obtain Chatzidakis' Theorem for the finite $\SU$-rank case, apart from the decomposition in orthogonal components:

\kor Let $b=\bdd(\cb(a/b))$ be such that $\SU(b)<\omega$. Then $\tp(b/\bdd(b)\cap\bdd(a))$ is  analysable in the family of all non one-based types of $\SU$-rank $1$.\ekor

\section{Applications and the Canonical Base Property}
In this section and the next, $\Sigma^{nob}$ will be the family of non one-based regular types (seen as partial types). 
For the applications one would like (and often has) more than mere strongly $\Sigma^{nob}$-basedness of canonical bases:
\begin{definition} A supersimple theory $T$ has the {\em Canonical Base Property CBP} if $\tp(\cb(a/b)/a)$ is almost $\Sigma^{nob}$-internal for all $a$, $b$.\end{definition}
\bem In other words, in view of Corollary \ref{sigmabased} a theory has the CBP if for all $a,b$ the type $\tp(\cb(a/b)/a)$ is $\Sigma^{nob}$-flat.\ebem
It had been conjectured that all supersimple theories of finite rank have the CBP, but Hrushovski has constructed a counter-example \cite{hpp}.
\bem Chatzidakis has shown for types of finite $SU$-rank that the CBP implies that even $\tp(\cb(a/b)/\bdd(a)\cap\bdd(b))$ is almost $\Sigma^{nob}$-internal \cite[Theorem 1.15]{zoe}.\ebem
\bsp The CBP holds for types of finite rank in
\begin{itemize}
\item Differentially closed fields in characteristic $0$ \cite{PZ03}.
\item Generic difference fields \cite{PZ03,zoe}.
\item Compact complex spaces \cite{camp,fuji,pi02}.\end{itemize}
Moreover, in those cases we have a good knowledge of the non one-based types.\ebsp

Kowalski and Pillay \cite[Section 4]{ko-pi} have given some consequences of strongly $\Sigma$-basedness in the context of 
groups. In fact, they work in a theory with the CBP, but they remark that their results hold, with {\em $\Sigma$-analysable} 
instead of {\em almost $\Sigma$-internal}, in all stable strongly $\Sigma$-based theories.
\tats\label{kopil} Let $G$ be an $\emptyset$-hyperdefinable strongly $\Sigma$-based group in
a stable theory.\begin{enumerate}
\item If $H\le G$ is connected with canonical parameter $c$, then $\tp(c)$ is $\Sigma$-analysable.
\item $G/Z(G)$ is $\Sigma$-analysable.
\end{enumerate}\etats

An inspection of their proof shows that mere simplicity of the ambient theory is sufficient, replacing centers by approximate centers and connectivity by local connectivity. Recall that the {\em approximate center} of a group $G$ is
$$\tilde Z(G)=\{g\in G:[G:C_G(g)]<\infty\}.$$
A subgroup $H\le G$ is {\em locally connected} if for all group-theoretic or model-theoretic conjugates $H^\sigma$ of $H$, if $H$ and $H^\sigma$ are commensurate, then $H=H^\sigma$. Locally connected subgroups and their cosets have canonical parameters; every subgroup is commensurable with a unique minimal locally connected subgroup, its {\em locally connected component}. For more details about the approximate notions, the reader is invited to consult \cite[Definition 4.4.9 and Proposition 4.4.10]{wa00}.
\satzli\label{kopilsimple} Let $G$ be an $\emptyset$-hyperdefinable strongly $\Sigma$-based
group in a
simple theory.\begin{enumerate}
\item If $H\le G$ is locally connected with canonical parameter $c$, then $\tp(c)$ is $\Sigma$-analysable.
\item $G/\tilde Z(G)$ is $\Sigma$-analysable.
\end{enumerate}\esatzli
\bew (1) Take $h\in H$ generic over $c$ and $g\in G$ generic over $c,h$. Let $d$
be the canonical parameter of $gH$. Then $\tp(gh/g,c)$ 
is the generic type of $gH$, so $d$ is interbounded with $\cb(gh/g,c)$. By
strongly $\Sigma$-basedness, $\tp(d/gh)$ is $\Sigma$-analysable. But
$c\in\dcl(d)$, so $\tp(c/gh)$ is $\Sigma$-analysable, as is $\tp(c)$ since
$c\ind gh$.

(2) For generic $g\in G$ put
$$H_g=\{(x,x^g)\in G\times G:x\in G\},$$
and let $H_g^{lc}$ be the locally connected component of $H_g$. Then $g\tilde Z(G)$ is interbounded with the canonical parameter of $H_g^{lc}$, so $\tp(g\tilde Z(G))$ is $\Sigma$-analysable, as is $G/\tilde Z(G)$.\qed

\satz\label{GoverN} Let $G$ be an $\emptyset$-hyperdefinable strongly $\Sigma$-based group in a simple theory. If $G$ is supersimple or type-definable, there is a normal nilpotent $\emptyset$-hyperdefinable subgroup $N$ such that $G/N$ is almost $\Sigma$-internal. In particular, a supersimple or type-definable group $G$ in a simple theory has a normal nilpotent hyperdefinable subgroup $N$ such that $G/N$ is almost $\Sigma^{nob}$-internal.\esatz
\bew $G/\tilde Z(G)$ is $\Sigma$-analysable by Proposition \ref{kopilsimple}. Hence there is a continuous sequence
$$G=G_0\triangleright G_1\triangleright G_2\triangleright\cdots\triangleright G_\alpha\triangleright\tilde Z(G)$$
of normal $\emptyset$-hyperdefinable subgroups such that successive quotients $Q_i=G_i/G_{i+1}$ are $\Sigma$-internal for all $i<\alpha$, and $G_\alpha/\tilde Z(G)$ is bounded.

Now $G$ acts on every quotient $Q_i$. Let
$$N_i=\{g\in G:[Q_i:C_{Q_i}(g)]<\infty\}$$
be the approximate stabilizer of $Q_i$ in $G$, again an $\emptyset$-hyperdefinable subgroup. If $(q_j:j<\kappa)$ is a long independent generic sequence in $Q_i$ and $g$, $g'$ are two elements of $G$ which have the same action on all $q_j$ for $j<\kappa$, there is some $j_0<\kappa$ with $q_{j_0}\ind g,g'$. Since $g^{-1}g'$ stabilizes $q_{j_0}$ it lies in $N_i$, and $gN_i$ is determined by the sequence $(q_j,q_j^g:j<\kappa)$. Thus $G/N_i$ is $Q_i$-internal, whence $\Sigma$-internal.

Put $N=\bigcap_{i<\alpha}N_i$. Since $\prod_{i<\alpha}G/N_i$ projects definably onto $G/N$, the latter quotient is also $\Sigma$-internal. In order to finish it now suffices to show that $N$ is virtually nilpotent. In particular, we may assume that $N$ is $\emptyset$-connected.

Consider the approximate ascending central series $\tilde Z_i(N)$. Note that $N$ centralizes $G_\alpha/\tilde Z(G)$ by $\emptyset$-connectivity. Moreover, $N$ approximatively stabilizes all quotients $(G_i\cap N)/(G_{i+1}\cap N)$. Hence, if $G_{i+1}\cap N\le\tilde Z_j(N)$, then $G_i\cap N\le\tilde Z_{j+1}(N)$.
If $G$ is supersimple, we may assume that all the $Q_i$ are unbounded, so $\alpha$ is finite and $N=\tilde Z_{\alpha+2}(N)$. In the type-definable case, note that $\tilde Z_i(N)$ is relatively $\emptyset$-definable by \cite[Lemma 4.2.6]{wa00}. So by compactness the least ordinal $\alpha_i$ with $G_{\alpha_i}\cap N\le\tilde Z_i(N)$ must be a successor ordinal, and $\alpha_{i+1}\le\alpha_i-1<\alpha_i$. Hence the sequence must stop and there is $k<\omega$ with $N=\tilde Z_k(N)$. But then $N$ is nilpotent by \cite[Proposition 4.4.10.3]{wa00}.\qed

\bem In a similar way one can show that if $G$ acts definably and faithfully on a $\Sigma$-analysable group $H$ and $H$ is supersimple or type-definable, then there is a hyperdefinable normal nilpotent subgroup $N\triangleleft G$ such that $G/N$ is almost $\Sigma$-internal.\ebem

\section{Final Remarks}
We have seen that for (weak) $\Sigma$-ampleness only the first level of an element is important. However, the difference between strong $\Sigma^{nob}$-basedness and the CBP is precisely the possible existence of a second (or higher) $\Sigma^{nob}$-level of $\cb(a/b)$ over $a$, i.e.\ its non $\Sigma^{nob}$-flatness.

One might be tempted to try to prove the CBP replacing $\Sigma^{nob}$-closure by $\Sigma^{nob}_1$-closure. In fact it is possible to define a corresponding notion of $\Sigma_1$-ampleness, and to prove an analogue of Theorem \ref{analysable}. However, since $\Sigma_1$-closure is not a closure operator, the equivalence between $\Sigma_1^{nob}$-basedness (i.e.\ the CBP) and non $1$-$\Sigma_1^{nob}$-ampleness breaks down. So far we have not found a way around this.

A possible approach to circumvent the failure of the CBP in general could be to use Theorem \ref{GoverN} in the applications, rather than establish the CBP for particular theories and use Fact \ref{kopil} (or Proposition \ref{kopilsimple}), but we have not looked into this.

Finally, it might be interesting to look for a variant of ampleness  which does take all levels into account, as one might hope to obtain stronger structural consequences.

\end{document}